\documentclass[12pt,twoside]{article}

\usepackage{epsfig} \epsfverbosetrue
\input epsf
\usepackage{headerfo}
\usepackage{amsfonts}

\newcounter{theorem}
\newcommand{\newsection}[1]{{\setcounter{theorem}{0} \section{#1}}}
\newtheorem{Theorem}{Theorem}[section]
\newtheorem{Definition}[Theorem]{Definition}
\newtheorem{Proposition}[Theorem]{Proposition}
\newtheorem{Lemma}[Theorem]{Lemma}

\def\eop{{ \vrule height7pt width7pt depth0pt}\par\bigskip}

\newcommand{\R}{\mathbb R}
\newcommand{\C}{\mathbb C}

\newif\ifpdf
\ifx\pdfoutput\undefined
  \pdffalse
\else
  \pdfoutput=1
  \pdftrue
\fi

\usepackage{latexsym}

\chardef\aa=64

\begin{document}

\renewcommand{\author}{T.~Bloom \\  
Department of Mathematics \\
University of Toronto\\
Toronto, Ontario\\
Canada M5S 2E4,\\
\medskip L.~Bos \\  
Department of Mathematics and Statistics\\
University of Calgary \\
Calgary, Alberta\\
Canada T2N 1N4, \\
\medskip
N.~Levenberg\\
Department of Mathematics \\
Indiana University \\
Bloomington, Indiana, \\
USA\\
\medskip
and \\
\medskip
S.~Waldron \\
Department of Mathematics\\
University of Auckland\\
Auckland, New Zealand}

\newcommand{\stitle}{Convergence of Optimal Measures}
 
\renewcommand{\title}{On the Convergence of Optimal Measures}

\renewcommand{\date}{\today}
\flushbottom
\setcounter{page}{1}
\pageheaderlinetrue
\oddpageheader{}{\stitle}{\thepage}
\evenpageheader{\thepage}{\stitle}{}
\thispagestyle{empty}
\vskip1cm
\begin{center}
\LARGE{\bf \title}
\\[0.7cm]
\large{\author}
\end{center}
\vspace{0.3cm}
\begin{center}
\date
\end{center}
\vfill \eject
\begin{abstract}
Using recent results of Berman and Boucksom [BB2] we show that for a non-pluripolar compact set $K\subset\C^d$ and an admissible weight function $w=e^{-\phi}$ any sequence of so-called optimal measures converges
weak-* to the equilibrium measure $\mu_{K,\phi}$ of (weighted)  Pluripotential Theory for $K,\phi$.
\end{abstract}

\setcounter{section}{0}

\newsection{Introduction} 
Suppose that $K\subset \C^d$ is compact and non-pluripolar (see Definition 2.1) and that $\mu$ is a probability measure on $K.$ We let  ${\cal P}_n$ denote the holomorphic polynomials of degree at most $n$ 
and we assume that $\mu$ is non-degenerate on ${\cal P}_n,$ i.e., with the associated inner-product
\begin{equation}\label{ip}
\langle f,g\rangle_\mu:=\int_K f\overline{g}d\mu 
\end{equation}
and $L_2(\mu)$ norm, $\|f\|_{L_2(\mu)}=\sqrt{\langle f,f\rangle_\mu}$, we have $\|p\|_{L_2(\mu)}=0$ for $p\in {\cal P}_n$ implies that $p=0$. It follows from the reasoning used in Proposition 3.5 of [B] that $\mu$ is
non-degenerate on ${\cal P}_n$ if and only if supp$(\mu)$ is not contained in an algebraic variety of degree $n$. 
Then ${\cal P}_n$ equipped with the inner-product (\ref{ip}) is a finite dimensional Hilbert
space of dimension  
\[ N:={d+n \choose n}.\]
We may also consider the uniform norm on $K,$
\[\|f\|_K:=\max_{z\in K}|f(z)|\]
and it is natural to compare the two norms for $p\in{\cal P}_n.$

Since $\mu$ is a probability measure we always have
\[\|p\|_{L_2(\mu)}\le \|p\|_K.\]
Moreover since ${\cal P}_n$ is finite dimensional there is always a constant $C=C(n,\mu,K)$ such that
the reverse inequality holds,
\[\|p\|_K\le C \|p\|_{L_2(\mu)}.\]
In fact, as is well known and easy to verify, the {\it best} constant $C$ (also sometimes called the
Berstein-Markov factor)
is given by 
\[C=\sup_{p\in {\cal P}_n,\,p\ne0}{\|p\|_K\over \|p\|_{L_2(\mu)}}=\max_{z\in K}\sqrt{K_n^\mu(z)}\]
where 
\begin{equation}\label{Kn}
K_n^\mu(z):=\sum_{j=1}^N|q_j(z)|^2
\end{equation}
is the diagonal of the reproducing kernel for ${\cal P}_n,$ sometimes also called the
(reciprocal of the) Christoffel function, and 
$Q_n=\{q_1,q_2,\cdots,q_N\}$ is an orthonormal basis for ${\cal P}_n.$ 

It is natural to ask among all probability measures on $K,$ which one provides the smallest such factor.
\begin{Definition}
Suppose that the probability measure $\mu$ has the property that 
\[\max_{z\in K}\sqrt{K_n^\mu(z)}\le \max_{z\in K}\sqrt{K_n^{\mu'}(z)}\]
for all other probability measures $\mu'$ on $K.$ Then we say that
$\mu$ is an optimal measure of degree $n$ for $K.$
\end{Definition}

Note that for {\it any} probabilty measure $\mu,$  $\displaystyle{\int_K K_n^\mu(z) d\mu=N},$ so that
\[\max _{z\in K} K_n^\mu(z)\ge N.\]
It turns out that for an Optimal Measure (see Lemma \ref{MaxIsN}  below)
\begin{equation}\label{max=N}
\max _{z\in K} K_n^\mu(z) = N.
\end{equation}

\subsection{A Second Optimality Property}

Optimal Measures also enjoy a second extremal property. To see this
let
\[B_n=\{p_1,p_2,\cdots, p_N\}\]
be a (fixed) basis for ${\cal P}_n$ and consider the associated Gram matrix
\begin{equation}\label{GramMatrix}
G_n^\mu(B_n):=[\langle p_i,p_j\rangle_\mu]\in\C^{N\times N}.
\end{equation}
If we expand $p_i$ in the orthonormal basis $Q_n$ we obtain
\begin{equation}\label{Pexpanded}
p_i=\sum_{k=1}^N\langle p_i,q_k\rangle_\mu q_k
\end{equation}
so that
\begin{eqnarray*}
\langle p_i,p_j\rangle_\mu&=&\sum_{k=1}^N \langle p_i,q_k\rangle_\mu \langle q_k,p_j\rangle_\mu \\
&=& \sum_{k=1}^N \langle p_i,q_k\rangle_\mu \overline{\langle p_j,q_k\rangle_\mu}.
\end{eqnarray*}
It follows that we have the factorization
\begin{equation}\label{factorG}
G_n^\mu(B_n)=V_n^\mu(V_n^\mu)^*
\end{equation}
where
\begin{equation}\label{Vn}
V_n^\mu=V_n^\mu(B_n,Q_n):=[\langle p_i,q_j\rangle_\mu]\in\C^{N\times N}.
\end{equation}
If now $\mu'$ is another probabilty measure on $K$ with associated inner-product 
$\langle f,g\rangle_{\mu'}$ and orthonormal basis $Q'_n=\{q_1',q_2',\cdots,q_N'\},$ then
from the expansion (\ref{Pexpanded}) we obtain
\begin{eqnarray*}
(V_n^{\mu'})_{ij}&=&\langle p_i,q_j'\rangle_{\mu'} \\
&=&\sum_{k=1}^N \langle p_i,q_k\rangle_\mu \langle q_k, q_j'\rangle_{\mu'} \\
&=& \sum_{k=1}^N (V_n^\mu)_{ik} A_{kj}
\end{eqnarray*}
where
\begin{equation}\label{TransitionMatrix}
A=A(Q_n,Q_n',\mu,\mu'):=[\langle q_i,q_j'\rangle_{\mu'}]\in\C^{N\times N}.
\end{equation}
Hence we have the transition
\begin{equation}\label{Trans}
V_n^{\mu'}=V_n^{\mu}A.
\end{equation}
Now, the transition matrix has the property that
\begin{eqnarray*}
\sum_{i=1}^N\sum_{j=1}^N|A_{ij}|^2 &=&\sum_{i=1}^N \left\{\sum_{j=1}^N |\langle q_i,q_j'\rangle_{\mu'}|^2\right\} \\
&=& \sum_{i=1}^N |\langle q_i,q_i\rangle_{\mu'}|^2\quad \hbox{(by Parseval)}\\
&=& \sum_{i=1}^N \int_K |q_i(z)|^2 d\mu' \\
&=& \int_K K_n^\mu(z) d\mu'.
\end{eqnarray*}
Hence if $\mu$ is an Optimal Measure, satisfying (\ref{max=N}), we have
\[ {\rm tr}(A^*A)=\sum_{i=1}^N\sum_{j=1}^N|A_{ij}|^2 \le N \]
for {\it any} other probability measure $\mu'.$ From this it follows that
the sum of the eigenvalues
\[\sum_{k=1}^N\lambda_k(A^*A)={\rm tr}(A^*A)\le N\]
and hence, by the Arithmetic-Geometric Mean inequality, 
\[{\rm det}(A^*A)=\prod_{k=1}^N\lambda_k(A^*A)\le \left({1\over N}\sum_{k=1}^N\lambda_k(A^*A)\right)^N\le 1,\]
i.e., if $\mu$ is an Optimal Measure satisfying (\ref{max=N}) and $\mu'$ is any other probability measure, then
the determinant of the transition matrix $A$ satisfies
\[|\det(A)|\le 1.\]
Consequently, by (\ref{Trans}), 
\[|{\rm det}(V_n^{\mu'})|\le |{\rm det}(V_n^\mu)|\]
and by the factorization (\ref{factorG})
\begin{equation}\label{OptimalDet}
|{\rm det}(G_n^{\mu'}(B_n))|\le |{\rm det}(G_n^\mu(B_n))|,
\end{equation}
i.e., an Optimal Measure $\mu$ also maximizes the determinant of the associated Gram matrix.

\subsection{Optimal Polynomial Interpolation}

There is a close connection between Optimal Measures and the so-called Fekete points of polynomial interpolation. Indeed,
suppose that $\mu$ is a discretely supported (probability) measure of the form
\begin{equation}
\label{equalwt}
\mu={1\over N}\sum_{i=1}^N\delta_{x_i},\quad x_i\in K.
\end{equation}
Then if $\mu$ is non-degenerate on the polynomials of degree $n,$ 
it is easy to see that $q_i=\sqrt{N}\ell_i,$ $1\le i\le N,$ where $\ell_i$ is the $i$th fundamental Lagrange
polynomial for the points $\{x_i\},$ form an orthonormal set with respect to $\langle \cdot,\cdot\rangle_\mu.$
Hence
\begin{eqnarray*}
(V_n^\mu)_{ij}&=&\langle p_i,q_j\rangle_\mu \\
&=& \sqrt{N} \langle p_i,\ell_j\rangle_\mu \\
&=& \sqrt{N} {1\over N}\sum_{k=1}^N p_i(x_k)\ell_j(x_k) \\
&=& {1\over\sqrt{N}} p_i(x_j)
\end{eqnarray*}
so that $V_n^\mu$ is in this case (a multiple of) the Vandermonde matrix for the basis $B_n$ and the points
$\{x_i\}.$ Hence maximizing $|{\rm det}(G_n^\mu)|$ over all discrete probability measures of the form (\ref{equalwt}) is equivalent to maximizing the modulus of the Vandermonde
determinant. Points which do this are called Fekete points for $K$ and the corresponding discrete measure is
said to be a Fekete measure.

With regard to the Christoffel function, we have 
\[K_n(z)=\sum_{k=1}^N|q_i(z)|^2=N\sum_{k=1}^N |\ell_k(z)|^2\]
so that minimizing $\max_{z\in K} K_n(z)$ over discrete measures of the form (\ref{equalwt}) is equivalent to finding the points 
for which $\max_{z\in K} \sum_{k=1}^N |\ell_k(z)|^2$ is as small as possible. This problem (for the interval
$K=[-1,1]$) was first studied by Fej\'er [F] and hence we refer to the solution points as
Fej\'er points and the corresponding measure as a Fej\'er measure. We remark that, in general, Fekete measures
and Fej\'er measures need not coincide (although they do in the univariate case of $K=[-1,1]$), cf. [Bo].

Further, if we regard the projection $\pi_\mu$ from $C(K)$ to ${\mathcal P}_n$ 
\[\pi_\mu(f):=\sum_{j=1}^N\langle f,q_j\rangle_\mu q_j=\sum_{j=1}^N f(x_j)\ell_j\]
as a map from $C(K)\to C(K),$ with both spaces equipped
with the {\it uniform} norm, then it is easy to see that
\[\|\pi_\mu\|=\Lambda_n:=\max_{z\in K}\sum_{k=1}^N|\ell_k(z)|,\]
the so-called Lebesgue constant for the interpolation process. The points for which $\Lambda_n$ is as small
as possible are called the Lebesgue points and will in general be different from both the Fekete and Fej\'er
points. We return to Lebesgue constants in a remark at the end of the paper.

\subsection{Optimal Experimental Designs}

Consider a polynomial $p\in {\cal P}_n$ which we write in the form
\[p=\sum_{k=1}^N \theta_k p_k.\]
Suppose that we observe the values of $p$ at $M\ge N$ points $x_j\in K$ with some random errors, i.e., we observe
\[y_j =p(x_j)+\epsilon_j,\quad 1\le j\le N\]
where we assume that the errors $\epsilon_j\sim N(0,\sigma)$ are independent. In matrix form this becomes
\[ y= X\theta+\epsilon\]
where $y,\theta,\epsilon\in\C^N$ and
\[X=\left[\begin{array}{cccccc}
p_1(x_1)&p_2(x_1)&\cdot&\cdot&\cdot&p_N(x_1) \cr
p_1(x_2)&p_2(x_2)&\cdot&\cdot&\cdot&p_N(x_2) \cr
\cdot&&&&&\cdot\cr
\cdot&&&&&\cdot\cr
\cdot&&&&&\cdot\cr
\cdot&&&&&\cdot\cr
\cdot&&&&&\cdot\cr
p_1(x_M)&p_2(x_M)&\cdot&\cdot&\cdot&p_N(x_M) \end{array}\right]\in \C^{M\times N}.\]
Our assumption on the error vector $\epsilon$ means that
\[{\rm cov}(\epsilon)=\sigma^2I_N\in\R^{N\times N}.\]
Now, the least squares estimate of $\theta$ is
\[\widehat{\theta}:=(X^*X)^{-1}X^*y\]
and we may compute the covariance matrix 
\[{\rm cov}(\widehat{\theta})=\sigma^2(X^*X)^{-1}.\]
Hence the confidence region of level $t$ for $\theta$ is the set
\begin{eqnarray*}
&&\{\theta\in \C^N\,:\,(\theta-\widehat{\theta})^*[{\rm cov}(\widehat{\theta})]^{-1}(\theta-\widehat{\theta})\le t\} \\
&=&\{\theta\in \C^N\,:\,\sigma^{-2}(\theta-\widehat{\theta})^*(X^*X)(\theta-\widehat{\theta})\le t\}.
\end{eqnarray*}
The volume of such a set is proportional to $1/\sqrt{{\rm det}(X^*X)}$ and hence maximizing the
${\rm det}(X^*X)$ is equivalent to choosing the observation points $x_i\in K$ so as to have the most
``concentrated'' confidence region for the parameter to be estimated.

Note however that the entries of $\displaystyle{ {1\over M}X^*X}$ are the discrete inner products of the $p_i$
with respect to the measure
\begin{equation}\label{StatsMeas}
\mu ={1\over M}\sum_{k=1}^M \delta_{x_k},
\end{equation}
i.e., $\displaystyle{ {1\over M}X^*X}$ is the Gram matrix associated to this $\mu.$ Hence we may think of an Optimal Measure as that which gives the confidence region of greatest concentration.

There is also a second statistical interpretation of Optimal Measures.
If we set
\begin{equation}\label{P}
P(x)=\left[\begin{array}{c}p_1(x)\cr p_2(x)\cr\cdot\cr\cdot\cr p_N(x)\end{array}\right]\in \C^{ N}
\end{equation}
then the least squares estimate of the observed polynomial is
\[P^t(x)\widehat{\theta}.\]
We may compute its variance to be
\begin{eqnarray*}
{\rm var}(P^t(x)\widehat{\theta})&=&\sigma^2P^*(x)(X^*X)^{-1}P(x) \\
&=&{1\over M}\sigma^2 P^*(x)(G_n^\mu)^{-1}P(x)
\end{eqnarray*}
where $\mu$ is again given by (\ref{StatsMeas}). But, it is not difficult to see that
\[P^*(x)(G_n^\mu)^{-1}P(x)=K_n^\mu(x)\] 
so that
\[{\rm var}(P^t(x)\widehat{\theta})={1\over M}\sigma^2K_n^\mu(x)\]
and the experiment that minimizes the maximum variance of the estimate of the observed polynomial is exactly the
one that minimizes the maximum of $K_n^\mu.$

We hope that the reader is convinced that Optimal Measures are interesting and worthy of further study. More about optimal experimental design may be found in the monographs [KS] and
[DS]. In the next section we introduce a slightly generalized (weighted) version of Optimal Measures and show that
they converge weak-* to the so-called equilibrium measure of Pluripotential Theory for $K.$

\newsection{Weighted Optimal Measures}
We recall the definition of plurisubharmonic function and pluripolar set. 

\begin{Definition}
A function $u\,:\,\C^d\to[-\infty,\infty)$ is said to be plurisubharmonic (psh) if it is upper semi-continuous (usc)
and, when restricted to any complex line, is either subharmonic or identically $-\infty.$ A set $E\subset \C^d$ is pluripolar if $E\subset \{z\in \C^d: u(z)=-\infty\}$ for some psh $u$ (with $u\not \equiv -\infty$).
\end{Definition}

Suppose that $K\subset \C^d$ is compact and non-pluripolar. 

\begin{Definition}
A function $w\,:\,K\to \R$ is said to be an admissible weight function if\par
\noindent (i) $w\ge 0$ on $K$\par
\noindent (ii) $w$ is upper semi-continuous\par
\noindent (iii) the set
\[\{z\in K\,:\, w(z)>0\}\]
is not pluripolar.
\end{Definition}
For $\mu$ a probability measure on $K$ and admissible weight $w$ we denote the associated weighted inner product
of degree $n$ by
\begin{equation}\label{Wip}
\langle f,g\rangle_{\mu,w}:=\int_K f(z)\overline{g(z)}w^{2n}(z)d\mu.
\end{equation}
For a (fixed) basis $B_n=\{p_1,p_2,\cdots, p_N\}$ of ${\cal P}_n$ we form the Gram matrix
\begin{equation}\label{WGramMatrix}
G_n^{\mu,w}=G_n^{\mu,w}(B_n):=[\langle p_i,p_j\rangle_{\mu,w}]\in\C^{N\times N}
\end{equation}
and the associated weighted Christoffel function
\begin{equation}\label{WKn}
K_n^{\mu,w}(z):=\sum_{j=1}^N|q_j(z)|^2w^{2n}(z)
\end{equation}
where, as before, $Q_n=\{q_1,q_2,\cdots,q_N\}$ is an orthonormal basis for ${\cal P}_n$ with respect
to the inner-product (\ref{Wip}). We note that as the Christoffel function is (essentially) the diagonal of the reproducing
kernel, it is independent of the particular orthonormal basis $Q_n.$

\begin{Definition}
Suppose that $w$ is an admissible weight on $K.$ If a probability measure $\mu$ has either of the two following, equivalent, properties:\par
(a) $\displaystyle{{\rm det}(G_n^{\mu',w})\le {\rm det}(G_n^{\mu,w})}$
for all other probability measures $\mu'$ on $K$\par
\noindent or\par
(b) $\displaystyle{\max_{z\in K}K_n^{\mu,w}(z)=N}$\par
\noindent then $\mu$ is said to be an Optimal Measure of
degree $n$ for $K$ and $w.$
\end{Definition}

By (the proof of) Lemma 2.1 of [KS, Chapter X], the set of matrices
\[\{G_n^{\mu,w}\,:\,\mu\,\,\hbox{is a probability measure on}\,\,K\}\]
is compact (and convex). Hence, by property (a), an Optimal Measure always exists. They need not be unique.

That the conditions (a) and (b) are equivalent (in the unweighted case) is the content
of the Kiefer-Wolfowitz Equivalence Theorem [KW] (but see also [KS, Theorem 2.1, Chapter X] or
else [Bo]).

Although the references cited prove this theorem only in the unweighted case, the generalization to the
weighted case is completely straightforward, and hence we do not include a separate proof. It is however
useful to note that (as is easy to see) with $P$ defined as in (\ref{P}),
\begin{equation}\label{OtherKn}
w^{2n}P^*(G_n^{\mu,w})^{-1}P = K_n^{\mu,w}.
\end{equation}

An important property of Optimal Measures is
\begin{Lemma}\label{MaxIsN}
Suppose that $\mu$ is optimal for $K$ and $w.$ Then
\[K_n^{\mu,w}(z)=N,\quad a.e.\,\, [\mu].\]
\end{Lemma}

\noindent {\bf Proof.} On the one hand
\[\max_{z\in K}K_n^{\mu,w}(z)=N\]
while on the other hand, by the orthonormality of the $q_j,$
\[ \int_K K_n^{\mu,w}(z)\, d\mu=\int_K\sum_{j=1}^N|q_j(z)|^2w^{2n}(z)\,d\mu(z)=N,\]
and the result follows. \eop

\medskip
We recall that for a basis $B_n$ and a set of points $Z_n=\{z_i\,:\,1\le i\le N\}\subset K$ the matrix
\[V_n=V_n(B_n,Z_n)=[p_i(z_j)]\in \C^{N\times N}\]
is called the Vandermonde matrix of the system. In case that the basis $B_n$ is the
{\it standard} monomial basis for ${\cal P}_n$ then we will write
\[ VDM(z_1,z_2,\cdots,z_N):={\rm det}(V_n).\]
Of fundamental importance for us will be
\begin{Definition}
Suppose that $K\subset \C^d$ is compact and that $w$ is an admissible weight function 
on $K.$ We set
\[\delta_n^w(K):=\left(\max_{z_i\in K}|VDM(z_1,\cdots,z_N)|w^n(z_1)w^n(z_2)\cdots w^n(z_N)\right)^{1/m_n}\]
where $m_n=dnN/(d+1)$ is the sum of the degrees of the $N$ monomials of degree at most $n.$
Then 
\[\delta^w(K)=\lim_{n\to\infty} \delta_n^w(K)\]
is called the Weighted Transfinite Diameter of $K.$ We refer to $\delta_n^w(K)$ as the $n$th order
weighted transfinite diameter of $K.$
\end{Definition}
A proof that this limit exists may be found in [BL] or [BB1]; it was first proved in the unweighted case ($w\equiv 1$; i.e., $\delta^1(K)$) by Zaharjuta [Z].

Given the close connection between Vandermonde matrices and Gram matrices, as explained in the Introduction,
it is perhaps not suprising that we have
\begin{Proposition}\label{tfd}
Suppose that $K$ is compact that $w$ is an admissible weight function. Suppose further that
$\mu_n$ is an Optimal Measure of degree $n$ for $K$ and $w.$ Take the basis $B_n$ to be the standard
basis of monomials for ${\cal P}_n.$ Then
\[\lim_{n\to\infty}{\rm det}(G_n^{\mu_n,w})^{1/(2m_n)}=\delta^w(K).\]
\end{Proposition}
\noindent {\bf Proof.}
We first note the formula (cf. formula (3.3) of [BL])
\begin{eqnarray}
&&\int_{K^N}|VDM(z_1,\cdots,z_N)|^2w(z_1)^{2n}\cdots w(z_N)^{2n} d\mu_n(z_1)\cdots d\mu_n(z_N)\nonumber \\
&&=N!\,{\rm det}(G_n^{\mu_n,w}). \label{VDMint}
\end{eqnarray}
It follows immediately, since $\mu_n$ is a probability measure, that
\begin{equation}\label{ub}
{\rm det}(G_n^{\mu_n,w})\le {1\over N!}(\delta_n^w(K))^{2m_n}.
\end{equation}
Secondly, note that if $f_1,f_2,\cdots,f_N\in K$ are so-called (weighted) Fekete points of degree $n$ for $K,$
i.e., points in $K$ for which 
\[|VDM(z_1,\cdots,z_N)|w^n(z_1)w^n(z_2)\cdots w^n(z_N)\]
is maximal, then the discrete measure 
\begin{equation}\label{WFekete}
\nu_n={1\over N}\sum_{k=1}^N \delta_{f_k}
\end{equation}
based on these points is a candidate probability measure for property (a) of Definition 2.3. Hence
\[ {\rm det}(G_n^{\nu_n,w})\le {\rm det}(G_n^{\mu_n,w}).\]
But, as is easy to see,
\begin{eqnarray*}
{\rm det}(G_n^{\nu_n,w})&=&{1\over N^N}|VDM(f_1,\cdots,f_N)|^2w(f_1)^{2n}w(f_2)^{2n}\cdots w(f_N)^{2n} \\
&=&\left(\max_{z_i\in K}|VDM(z_1,\cdots,z_N)|w^n(z_1)w^n(z_2)\cdots w^n(z_N)\right)^2 \\
&=&{1\over N^N}(\delta_n^w(K))^{2m_n}.
\end{eqnarray*}
Hence, 
\[{1\over N^N}(\delta_n^w(K))^{2m_n} \le {\rm det}(G_n^{\mu_n,w}) \]
and combining this lower bound with the upper bound (\ref{ub}) we obtain
\[ {1\over N^N}(\delta_n^w(K))^{2m_n} \le {\rm det}(G_n^{\mu_n,w}) \le {1\over N!}(\delta_n^w(K))^{2m_n}\]
and the result follows. \eop

\medskip
Of course, it then follows that
\[\lim_{n\to \infty}{1\over 2 m_n}\log\,{\rm det}(G_n^{\mu_n,w})=\log(\delta^w(K)).\]
Now, suppose that $u\in C(K)$ and that $w(z)$ is an admissible weight function. Consider the weight $w_t(z):=w(z)\exp(-tu(z)),$ $t\in\R,$ and let $\mu_n$ be an optimal measure of degree $n$ for $K$ and $w.$
We set
\begin{equation}\label{fn}
f_n(t):=-{1\over 2 m_n}\log\,{\rm det}(G_n^{\mu_n,w_t}).\end{equation}
For $t=0,$ $w_0=w$ and hence we have
\[\lim_{n\to\infty}f_n(0)=-\log(\delta^{w}(K)).\]

\begin{Lemma}\label{1stderiv}We have
\[ f_n'(t)={d+1\over dN}\int_K u(z)K_n^{\mu_n,w_t}(z)d\mu_n.\]
In particular,
\begin{eqnarray}
f_n'(0)&=&{d+1\over dN}\int_K u(z)K_n^{\mu_n,w}(z)d\mu_n \nonumber \\
&=& {d+1\over d}\int_K u(z)d\mu_n \quad \hbox{(by Lemma \ref{MaxIsN})}.
\end{eqnarray}
\end{Lemma}
\noindent {\bf Proof.}
We calculate
\begin{eqnarray*}
2m_nf_n'(t)&=&-{d\over dt}{\rm trace}\left(\log(G_n^{\mu_n,w_t})\right) \\
&=& -{\rm trace} \left({d\over dt}\log(G_n^{\mu_n,w_t})\right) \\
&=& -{\rm trace}\left( (G_n^{\mu_n,w_t})^{-1}{d\over dt} G_n^{\mu_n,w_t}\right) \\
&=&2n\, {\rm trace}\left( (G_n^{\mu_n,w_t})^{-1} \left[\int_K p_i(z)\overline{p_j(z)}u(z)w(z)^{2n}\exp(-2ntu(z))d\mu_n\right]\right) \\
&=&2n\int_K P^*(z)(G_n^{\mu_n,w_t})^{-1}P(z) u(z)w(z)^{2n}\exp(-2ntu(z))d\mu_n\\
&=&2n\int_K u(z) P^*(z)(G_n^{\mu_n,w_t})^{-1}P(z) w_t(z)^{2n} d\mu_n \\
&=& 2n\int_K u(z) K_n^{\mu_n,w_t}(z) d\mu_n
\end{eqnarray*}
where the last equality follows from the remark (\ref{OtherKn}).

The result follows from the fact that $m_n=dnN/(d+1).$ \eop

\begin{Lemma}\label{2ndderiv}
The functions $f_n(t)$ are concave, i.e., $f_n''(t)\le0.$
\end{Lemma}
\noindent {\bf Proof.} First, let
\[g_n(h):=2m_nf_n(t+h)\]
so that $\displaystyle{f_n''(t)={1\over 2m_n}g_n''(0).}$ Also, note that if we change the
basis $B_n=\{p_1,\cdots,p_N\}$ to $C_n:=\{q_1,\cdots,q_N\}$ by $p_i=\sum_{j=1}^Na_{ij}q_j,$ 
then the Gram matrices transform (see e.g. [D, \S8.7]) by
\[G_n^{\mu_n,w_t}(B_n)=AG_n^{\mu_n,w_t}(C_n)A^*\]
where $A=[a_{ij}]\in \C^{N\times N}.$ Hence,
\[ g_n(h)=-\log({\rm det}(G_n^{\mu_n,w_{t+h}}(B_n)))=-\log({\rm det}(G_n^{\mu_n,w_{t+h}}(C_n)))-\log(|{\rm det}(A)|^2)\]
and we see that the derivatives of $g_n$ are independent of the basis chosen.

Let us choose $C_n$ to be an orthonormal basis for ${\cal P}_n$ with respect to the
inner-product $\langle \cdot,\cdot\rangle_{\mu_n,w}=\langle \cdot,\cdot\rangle_{\mu_n,w_0}.$ 

Now, for convenience, write $G(h)=G_n^{\mu_n,w_{t+h}}$ and set $F(h)=\log(G(h))$ so
that $G(h)=\exp(F(h)).$ By our choice of basis $C_n$ we have $G(0)=I\in \C^{N\times N},$ 
the identity matrix, and
$F(0)=[0]\in\C^{N\times N},$ the zero matrix.
Then, (see e.g. [Bh, p. 311]),
\[{dG\over dh}={d\over dh}\exp(F(h))=\int_0^1e^{(1-s)F(h)}{dF\over dh}e^{sF(h)}ds.\]
In particular
\[{dG\over dh}(0)={dF\over dh}(0).\]
Further, 
\begin{eqnarray*}
{d^2G\over dh^2}&=&\int_0^1 \left\{\left[{d\over dh}e^{(1-s)F(h)}\right]{dF\over dh}e^{sF(h)}
+e^{(1-s)F(h)}{d^2F\over dh^2}e^{sF(h)} \right. \\
&& \left. \quad +e^{(1-s)F(h)}{dF\over dh}\left[{d\over dh}e^{sF(h)}\right] \right\}ds.
\end{eqnarray*}
Evaluating at $h=0,$ using the fact that $F(0)=[0],$ we obtain
\begin{eqnarray*}
{d^2G\over dh^2}(0)&=&\int_0^1 \left\{(1-s){dF\over dh}(0)\times {dF\over dh}(0)\times I
+I\times {d^2F\over dh^2}(0)\times I \right.\\
&&\left. \quad+I\times {dF\over dh}(0)\times s{dF\over dh}(0)\, \right\}ds\\
&=&\int_0^1 \left\{(1-s+s)\left({dF\over dh}(0)\right)^2+{d^2F\over dh^2}(0)\right\}\,ds \\
&=&\left({dF\over dh}(0)\right)^2+{d^2F\over dh^2}(0).
\end{eqnarray*}
Hence,
\begin{eqnarray*}
{d^2F\over dh^2}(0)&=&{d^2G\over dh^2}(0)-\left({dF\over dh}(0)\right)^2\\
&=&[\int_Kq_i(z)\overline{q_j(z)}(-2nu(z))^2w_t(z)^{2n}d\mu_n]-[\int_Kq_i(z)\overline{q_j(z)}(-2nu(z))w_t(z)^{2n}d\mu_n]^2.
\end{eqnarray*}
It follows that
\begin{eqnarray*}
g_n''(0)&=&-{\rm trace}\left([\int_Kq_i(z)\overline{q_j(z)}(-2nu(z))^2w_t(z)^{2n}d\mu_n]\right)\\
&&+{\rm trace}\left([\int_Kq_i(z)\overline{q_j(z)}(-2nu(z))w_t(z)^{2n}d\mu_n]^2\right)\\
&=&-\sum_{i=1}^N\int_K|q_i(z)|^2w_t(z)^{2n}(2nu(z))^2d\mu_n \\
&&+\sum_{i=1}^N\sum_{j=1}^N\left|\int_Kq_i(z)\overline{q_j(z)}w_t(z)^{2n}(2nu(z))d\mu_n\right|^2\\
&=&-\sum_{i=1}^N \left\{\int_K|q_i(z)|^2w_t(z)^{2n}(2nu(z))^2d\mu_n - \right. \\
&&\quad \left.\sum_{j=1}^N \left|\int_Kq_i(z)\overline{q_j(z)}w_t(z)^{2n}(2nu(z))d\mu_n\right|^2\right\}.
\end{eqnarray*}
But notice that $\displaystyle{\int_Kq_i(z)\overline{q_j(z)}w_t(z)^{2n}(2nu(z))d\mu_n}$ is the
$j$th Fourier coefficient of the function $2nu(z)q_i(z)$ with respect to the orthonormal basis
$C_n,$ and also that $\displaystyle{\int_K|q_i(z)|^2w_t(z)^{2n}(2nu(z))^2d\mu_n}$ is the $L_2$ norm
squared of this same function. Hence, by Parsevals inequality,
\[g_n''(0)\le 0.\]
\eop
\newsection{The Limit of Optimal Measures}
Suppose again that $K\subset\C^d$ is compact and that $w$ is an admissible weight function. We set
\[\phi:=-\log(w).\]
In order
to state the convergence theorem we will need to briefly review some notions from (weighted) Pluripotential Theory.
We refer the reader to the monograph [K] and also to Appendix B of [ST] for more details.

The class of psh functions of at most logarithmic growth at infinity is denoted by
\[ {\cal L}:=\{u\,:\, u\,\,\hbox{is psh and}\,\,u(z)\le \log^+|z|+C\}.\]

Of special importance is the weighted {\it pluricomplex Green's function} (also known as the weighted
extremal function),
\begin{equation}\label{V}
V_{K,\phi}(z):=
\sup\,\{u(z)\,:\, u\in{\cal L},\,\,u\le \phi\,\,\hbox{on}\,\,K\}.
\end{equation}
The function $V_{K,\phi}^*(z)$ denotes the usc regularization of $V_{K,\phi}.$

Associated to the extremal function is the so-called {\it weighted equilibrium measure},
\begin{equation}\label{wem}
\mu_{K,\phi}:={1\over (2\pi)^d}(dd^cV_{K,\phi}^*)^d.
\end{equation}
Here $(dd^cv)^d$ refers to the Monge-Ampere operator (applied to $v$). That $\mu_{K,\phi}$ exists and is a probability measure
is one of the main results of Pluripotential Theory; we again refer the reader to [K] or Appendix B of [ST] for the details. We simply write $\mu_K$ in the unweighted case, i.e., $w=1$ and $\phi=0$.  
We remark, that in one variable, for $K=[-1,1]\subset\C,$ 
\[ \mu_{K}={1 \over \pi}{1\over\sqrt{1-x^2}}dx.\]

In a remarkable sequence of papers Berman and Boucksom [Be, BB1, BB2] have recently shown that the discrete
probability measures based on the weighted Fekete points (\ref{WFekete}) tend weak$-*$ to $\mu_{K,\phi}.$ Indeed, if for each $n$, $x_1^{(n)},x_2^{(n)},\cdots,x_N^{(n)}\in K$ are points in $K$ for which 
\[ \lim_{n\to \infty}\bigl[|VDM(x_1^{(n)},\cdots,x_N^{(n)})|w(x_1^{(n)})^nw(x_2^{(n)})^n\cdots w(x_N^{(n)})^n\bigr]^{1/m_n}=\delta^w(K)\]
({\it asymptotically} weighted Fekete points), then the discrete measures 
\[ \nu_n={1\over N}\sum_{k=1}^N \delta_{x_k^{(n)}}\]
converge weak$-*$ to $\mu_{K,\phi}.$ The main point of this note is to remark
that their proof may be extended to also give the limit of Optimal Measures. For completeness we give
the details of the proof, but we emphasize that it is their same argument as for the Fekete measure case.

\begin{Theorem}
Suppose that $K\subset \C^d$ is compact and that $w$ is an admissible weight function. We again set
$\phi:=-\log(w).$ Suppose further that $\mu_n$ is an Optimal Measure of degree $n$ for $K$ and $w.$
Then
\[\lim_{n\to\infty}\mu_n=\mu_{K,\phi}\]
where the limit is in the weak$-*$ sense.
\end{Theorem}
\noindent {\bf Proof.}
For $u\in C(K)$ we again set $w_t(z):=w(z)\exp(-tu(z))$ which corresponds to $\phi_t:=\phi +tu$ and $f_n(t)$ as in (\ref{fn}).
As mentioned above,
\[\lim_{n\to\infty} f_n(0)=-\log(\delta^w(K)).\]
Of fundamental importance is the Rumely formula for the transfinite diameter ([R, BB1]):
\begin{equation}\label{RumelyFormula}
-\log(\delta^w(K))={1\over d (2\pi)^d}{\cal E}(V_{K,\phi}^*,V_T).
\end{equation}
Here $V_T$ is the (unweighted) extremal function for a polydisc that contains $K$ and ${\cal E}$ is a
certain ``mixed energy'' whose exact formula is not important here. What is important is the
derivative formula of Berman and Boucksom [BB1],
\begin{equation}\label{BBderiv1}
\left.{d\over dt}{\cal E}(V_{K,\phi+tu},V_T)\right|_{t=0}=(d+1)\int_Ku(dd^cV_{K,\phi}^*)^d.
\end{equation}
In other words, setting $g(t)=-\log(\delta^{w_t}(K)),$
\begin{equation}\label{BBderiv2}
g'(0)={d+1\over d(2\pi)^d}\int_Ku(z)(dd^cV_{K,\phi}^*)^d.
\end{equation}
Now note that for each fixed $t,$ the measure $\mu_n,$ being optimal for $K$ and $w=w_0,$ is a
candidate for the optimal measure for $K$ and $w_t.$ If follows from property (a) of Definition 2.3 that
\[{\rm det}(G_n^{\mu_n,w_t}) \le {\rm det}(G_n^{\mu_n^t,w_t})\]
where we denote  an optimal measure for $K$ and $w_t$ by $\mu_n^t.$ Hence (see (\ref{fn})
\[ f_n(t)\ge -{1\over 2m_n}\log({\rm det}(G_n^{\mu_n^t,w_t}))\]
and consequently that
\begin{equation}\label{liminf}
\liminf_{n\to\infty} f_n(t) \ge -\log(\delta^{w_t}(K)) =g(t).
\end{equation}
It now follows from the simple Lemma \ref{BBLemma} below given by Berman and Boucksom [BB2] that
\[\lim_{n\to\infty}f_n'(0)=g'(0).\]
In other words, by Lemma \ref{1stderiv}, 
\begin{eqnarray*}
\lim_{n\to \infty}{d+1\over d}\int_K u(z)d\mu_n &=& {d+1\over d(2\pi)^d}\int_Ku(z)(dd^c V_{K,\phi}^*)^d \\
&=& {d+1\over d} \int_K u(z) d\mu_{K,\phi}.
\end{eqnarray*}
\eop
\begin{Lemma}\label{BBLemma}(Berman and Boucksom [BB2])
Let $f_n(t)$ be a sequence of concave functions on $\R$ and $g(t)$ a function on $\R.$ Suppose that
\[\liminf_{n\to\infty}f_n(t)\ge g(t),\quad \forall t\in\R\] and that
\[\lim_{n\to\infty} f_n(0)=g(0).\]
Suppose further that the $f_n$ and $g$ are differentiable at $t=0.$ Then
\[\lim_{n\to\infty}f_n'(0)=g'(0).\]
\end{Lemma}

\noindent{\bf Remark}. There exist many other natural sequences of measures $\{\mu_n\}$ which converge weak-* to $\mu_{K,\phi}$. For simplicity, we discuss the unweighted case ($\phi = 0$). Recall from subsection 1.2 that if $x_1,...,x_N\in K$, then $\Lambda_n:=\max_{z\in K}\sum_{k=1}^N|\ell_k(z)|$ 
is the so-called Lebesgue constant associated to polynomial interpolation at these points. Suppose for each $n=1,2,...$ we have $N=N(n)$ points $x_1^{(n)},...,x_N^{(n)} \in K$ with Lebesgue constant $\Lambda_n$. An elementary argument in [BBCL] shows that 
if $\limsup_{n\to \infty}\Lambda_n^{1/n} \leq 1$, then 
$$\lim_{n\to \infty} |VDM(x_1^{(n)},...,x_N^{(n)})|^{1/m_n}=\delta^1(K).$$
By the main result of [BB2], it follows that the discrete measures $$\mu_n :={1\over N}\sum_{i=1}^N\delta_{x_i^{(n)}}$$
converge weak-* to $\mu_K $. Since the fundamental Lagrange polynomials $\ell_i$ for Fekete points satisfy $||\ell_i||_K=1$ it is easy to see that the Lebesgue constants for either the Lebesgue 
or Fejer points satisfy the above growth condition so the weak-*  convergence to the equilibrium measure holds. Furthermore, in Proposition 3.7 of [BBCL] it was shown that for a so-called {\it Leja sequence} $\{x_1,x_2,...\}\subset K$, 
$$\lim_{n\to \infty} |VDM(x_1,...,x_N)|^{1/m_n}=\delta^1(K).$$
Thus, again from  [BB2] it follows that the discrete measures $$\mu_n :={1\over N}\sum_{i=1}^N\delta_{x_i}$$
converge weak-* to $\mu_K$. Such a sequence is defined inductively as follows. Take the standard monomial basis $\{p_1,p_2,...\}$ for $\cup_{n=0}^{\infty} \mathcal P_n$ ordered so that deg$p_i\leq$deg$p_j$ if $i\leq j$. Given $m$ points $z_1,...,z_m$ in $\C^d$, as before we write
$$VDM(z_1,...,z_m)=\det [p_i(z_j)]_{i,j=1,...,m}.$$
Starting with any point $x_1\in K$, having chosen $x_1,...,x_m\in K$ we choose $x_{m+1}\in K$ so that
$$|VDM(x_1,...,x_m,x_{m+1})|=\max_{x\in K} |VDM(x_1,...,x_m,x)|.$$
It is unknown if $\limsup_{n\to \infty}\Lambda_n^{1/n} \leq 1$ always holds for a Leja sequence, even if $d=1$.

\bigskip


\centerline {\bf REFERENCES}
\vskip10pt
\begin{description}

\item{[Be]} Berman, R., {\it Bergman Kernels for Weighted Polynomials and Weighted Equilibrium Measures
of $\C^n$}, preprint arXiv:math/0702357.

\item {[BB1]} Berman, R. and Boucksom, S., {\it Capacities and Weighted Volumes of Line Bundles}, preprint 
arXiv:0803.1950.

\item {[BB2]} Berman, R. and Boucksom, S., {\it Equidistribution of Fekete Points on Complex Manifolds},
preprint arXiv:0807.0035.

\item {[Bh]} Bhatia, R., {\bf Matrix Analysis}, GTM 169, Springer, New York, 1997.

\item{[B]} Bloom, T., {\it Oorthogonal polynomials in $\C^n$}, Indiana Univ. Math. J., Vol. 46, No. 2 (1997), 427 -- 452.

\item{[BBCL]} Bloom, T., Bos, L., Christensen, C. and Levenberg, N.,  {\it Polynomial interpolation of holomorphic functions  in $\C$ and $\C^n$}, Rocky Mtn. J. Math., 22 (1992),  441 -- 470.

\item{[BL]} Bloom, T. and Levenberg, N., {\it Transfinite diameter notions in $\C^N$ and integrals of Vandermonde
determinants}, preprint arXiv:0712.2844.

\item {[Bo]} Bos, L., {\it Some Remarks on the Fej\'er Problem for Lagrange Interpolation in Several Variables},
J. Approx. Theory, Vol. 60, No. 2 (1990), 133 -- 140.

\item{[D]} Davis, P.J., {\bf Interpolation and Approximation}, Dover, 1975.

\item{[DS]} Dette, H. and Studden, W.J., {\bf The Theory of Canonical Moments with Applications in Statistics,
Probability and Analysis}, Wiley Interscience, New York, 1997.

\item {[F]} Fej\'er, L., {\it Bestimmung dergenigen Abszissen eines Intervalles f\"ur welche die
Quadratsumme der Grundfunktionen der Lagrangeschen Interpolation im Intervalle eing moglichst kleines Maximum besitzt},
Ann. Scuoal Norm. Sup. Pisa (2) {\bf 1} (1932), 263 -- 276.

\item{[KS]} Karlin, S. and Studden, W.J.,  {\bf Tchebycheff Systems: With Applications in Analysis and Statistics}, Wiley Interscience, New York, 1966.

\item{[KW]} Kiefer, J. and Wolfowitz, J., {\it The equivalence of two extremum problems}, Canad. J. Math. {\bf 12}
(1960), 363 -- 366.

\item{[K]} Klimek, M., {\bf Pluripotential Theory}, Oxford Univ. Press, 1991.

\item{[R]} Rumely, R., {\it A Robin Foumula for the Fekete-Leja Transfinite Diameter}, Math. Ann. {\bf 337} no. 4 (2007), 729 -- 738.

\item{[ST]} Saff, E. and Totik, V., {\bf Logarithmic Potentials with External Fields}, Springer, 1997.

\item{[Z]} V. P. Zaharjuta, {\it Transfinite diameter, Chebyshev constants, and capacity for compacta in $\C^n$}, Math. USSR Sbornik, {\bf 25} (1975), no. 3, 350 -- 364.

\end{description}

\end{document}